\documentclass{amsart}

\usepackage{amsmath, amssymb}
\usepackage[utf8]{inputenc}
\usepackage{xcolor}
\usepackage{tikz}
\usepackage{tikz-cd}

\definecolor{dblue}{rgb}{0.0, 0.0, 0.55}

\usepackage{datetime2}

\usepackage[colorlinks, urlcolor=black]{hyperref}
\hypersetup{citecolor=dblue, linkcolor=dblue}

\numberwithin{equation}{section}

\title[On Kim's Assumption A]{On Kim's Assumption A over function fields}

\thanks{2020 \emph{Mathematics Subject Classification}. Primary 11F70, 22E50, 22E55.}

\author{H\'ector del Castillo}

\author{Luis Lomel\'i}

\date{\today}

\begin{document}

\maketitle

\parindent=17pt
\parskip=7pt

\begin{abstract}
We prove Kim's Assumtion A for the split classical groups in positive characteristic. Actually, we work in the slightly more general setting of groups of classical kind, which includes quasi-split classical groups and generalized spinor groups. We establish our results whenever a local Ramanujan bound holds; a bound that is known for the split classical groups in characteristic $p$, and we prove it for groups of classical kind under a local-global restriction.
\end{abstract}

\section{Introduction}

We work on the conjecture called ``Kim's Assumption A'', and settle the characteristic $p$ case for the split classical groups; and, in slightly more generality under a local-global restriction. The original assumption was introduced by Kim over number fields in \cite{Kim1999}, and further progress can be found in \cite{Kim2005,KimKim2011}. Given a globally generic cuspidal representation over a global field, it states that, at every place, the local normalized intertwining operator is holomorphic and non-zero on the complex right half plane starting at $1/2$.

The proof of Kim's Assumption A relies on three local ingredients for generic\footnote{We assume that generic representations are irreducible.} representations: the tempered $L$-function conjecture; the standard module conjecture; and, a Ramanujan bound discussed in \S~\ref{ss:class:Ram}. The first two ingredients are established for quasi-split groups over a non-archimedean local field of characteristic $p$ in \cite{dCHeLo}. The third ingredient, namely the local Ramanujan bound, is known for the split classical groups in characteristic $p$ \cite{Lom2019}. However, we work in slightly more generality by introducing groups of classical kind. These groups include the quasi-split classical groups and the generalized spinor groups. We state our results for groups of classical kind, while our proofs are made assuming the validity of the local Ramanujan bound; our results being unconditional for the split classical groups. 

An additional observation, for groups of classical kind, is that whenever the $\mathfrak{p}$-adic representation is a local component of a globally generic cuspidal representation, the local Ramanujan bound is also true. In characteristic zero, this was studied by Kim in \cite{Kim2000}. In positive characteristic, the case of quasi-split non-split orthogonal groups is studied in \cite{dC}. Here, we present the general result over function fields in \ref{ss:Kim:local-global}~Proposition.

We first state Kim's Assumption A in its original (global) form \ref{ss:KimA}, and then we follow it by a local variant \ref{ss:KimA:loc}; the local version implies the global one. The setting is that of the Langlands-Shahidi method, originally developed by Shahidi in characteristic zero, and extended to characteristic $p$ by the second named author. In this setting we have a pair of quasi-split connected reductive groups $\bf G$ and $\bf M$ defined over a global or a local field, where $\bf M$ is the Levi component of a proper maximal parabolic subgroup ${\bf P} = {\bf M}{\bf N}$ of $\bf G$ containing a fixed borel subgroup $\bf B$. The corresponding $L$-group  ${}^LM$, acts on the Lie algebra ${}^L\mathfrak{n}$ of ${}^L\!N$ by conjugation, and the irreducible constituents of the resulting representation are what we call LS representations; we label them by $r_i$, where the indexes $i$ vary from 1 to some $m$, according to the nilpotency class in ${}^L\mathfrak{n}$. 

Given a non-archimedean local field $F$ and a unitary generic representation $\pi$ of $M = \operatorname{\textbf M}(F)$, we have a local intertwining operator ${\rm A}(s,\pi)$ defined on the induced representation
\[ {\rm I}(s,\pi) = {\rm Ind}(\pi_{s\widetilde{\alpha}}), \]
where $\pi_{s\widetilde{\alpha}}$ is the representation $\pi$ after twisting by an unramified character of $M$ parametrized by a complex variable $s$; the notation is made precise in \S~\ref{s:prelim}.

Now, given a non-trivial character $\psi : F \rightarrow \mathbb{C}^\times$, we then have the normalizing factor
\[R(s,\pi,\psi)=\prod_{i=1}^{m} \frac{L(is,\pi,r_i)}{L(1+is,\pi,r_i)\varepsilon(is,\pi,r_i,\psi)}. \]
The normalized intertwining operator ${\rm N}(s,\pi)$ is defined to be such that
\begin{equation}\label{eq:N:def} A(s,\pi)=R(s,\pi,\psi){\rm N}(s,\pi)
\end{equation}
as a rational operator in $q^{-s}$. We now provide the statement of the assumption over global fields.

\subsection{Kim's Assumption A}\label{ss:KimA}
\emph{Let $K$ be a global field. Assume that $\bigotimes_{v} \pi_v$ is a globally generic unitary cuspidal automorphic representation of $\operatorname{\textbf M}(\mathbb A_K)$. Then ${\rm N}(s,\pi_v)$ is holomorphic and non-zero on $\operatorname{Re}(s)\geq 1/2$ for all places $v$.}

\medskip

We also find it useful to state a purely local version as follows.

\subsection{Local Kim's Assumption A}\label{ss:KimA:loc}
\emph{Let $F$ be a local field. Given a generic unitary representation $\pi$ of $M = {\bf M}(F)$, the normalized intertwining operator ${\rm N}(s,\pi)$ is holomorphic and non-zero on ${\rm Re}(s) \geq 1/2$.}

We observe that the latter version implies the former.

The present article begins with a section on preliminaries, where we explain the notation and recall the first two local ingredients that appear in the proof of Kim's Assumption A. The first is Shahidi's tempered $L$-function conjecture, which we state in \ref{ss:temp:L}. In characteristic zero, it was proved in most cases by Kim and Kim \cite{KimKim2011}, and was settled by Heiermann and Opdam in all cases \cite{HeOp2013}. Now in characteristic $p$, the authors establish the tempered $L$-function conjecture in collaboration with G. Henniart in \cite{dCHeLo}. The second local ingredient is the standard module conjecture, stated in \ref{ss:SM}. It was formulated by Casselman and Shahidi in \cite{CaSh1998}, who proved it for the classical groups. In general, over non-archimedean fields of characteristic zero, Heiermann and Mui\'c \cite{HeMu2007} proved the standard module conjecture assuming the validity of the tempered $L$-function conjecture, before it was proved to hold true. The characteristic $p$ case is now available, and can also be found in \cite{dCHeLo}.

The last ingredient is a local ``Ramanujan bound'' for the character appearing in the standard module. To specify and find this bound, we work with a group $\bf G$ of classical kind, a notion that we define in \S~\ref{ss:class:Ram}. Given a generic irreducible representation $\pi \simeq {\rm I}(\chi,\sigma)$ of $G$, expressed as a standard module, the character $\chi$ corresponds to a strictly increasing sequence of positive real numbers, known as the Langlans parameters of $\pi$. Then, we say that the Ramanujan bound is satisfied if these parameters are strictly less than $1$. In characteristic zero, this bound was studied by Kim in \cite{Kim2000}. Furthermore, for quasi-split classical groups the bound is a consequence of the Lapid-Mui\'c-Tadi\'c classification of generic unitary representations \cite{LaMuTa2004}. The Ramanujan bound for the split classical groups is transferred from characteristic $0$ to characteristic $p$ in \cite{Lom2019}, relying on the Kazhdan transfer of Ganapaty \cite{Ga2015}.

In our context of groups of classical kind, we prove the Ramanujan bound to be true whenever the representation is a local component of a globally generic cuspidal representation over a global function field, \ref{ss:Kim:local-global}~Proposition. In the course of proof, we encounter three different cases of LS representations and $L$-functions that merit special attention: Rankin-Selberg products of a pair of general linear group representations; a pair consisting of a representation of a general linear group and a representation of a group of classical kind; and finally, the case of a representation of a Siegel Levi subgroup. In each of these cases, we can apply the globalization of supercuspidal representations of general linear groups and Rankin-Selberg $L$-functions, and obtain the Ramanujan bound.

We arrive at \S~\ref{sec:kimproof} with all the ingredients ready to initialize the proof of Kim's Assumption A for characteristic $p$ local fields, the Ramanujan bound being the starting point. Over number fields, the methods were introduced by Kim, and we are guided over function fields by his arguments given in \cite{Kim1999, Kim2000, Kim2005,KimKim2011}. To prove Kim's Assumption A, we first reduce to the three cases mentioned above, by applying the standard module conjecture and the multiplicativity of Normalized intertwining operators. We then use the tempered $L$-function conjecture to prove the holomorphicity property. Finally, using Zhang's multiplicativity formula \cite{Zh}, we prove the non-vanishing property.

\subsection*{Acknowledgments} We would like to thank Guy Henniart, for enlightening mathematical discussions. The first author was supported by ANID Postdoctoral Project 3220656.  The second author was supported in part by FONDECYT Grant 1212013.

\section{Preliminaries}\label{s:prelim}

We let $F$ denote a non-archimedean local field throughout the article. We work with a quasi-split reductive group $\bf G$ defined over $F$, and a maximal parabolic subgroup ${\bf P} = {\bf M}{\bf N}$, with Levi component $\bf M$ and unipotent radical $\bf N$. We refer to such a Levi component simply as a maximal Levi of $\bf G$. Furthermore, we fix a borel subgroup $\bf B$ of $\bf G$ and we assume our parabolics contain $\bf B$, i.e., we work with standard parabolics. If $Z_{\bf G}$ is the center of $\bf G$, then $A_{\bf G}$ denotes its maximal split torus. Similarly, $A_{\bf M}$ denotes the corresponding maximal split torus of $\bf M$ inside $Z_{\bf M}$. Given an algebraic group $\bf H$, we use $H$ to denote the group of $F$ rational points of $\bf H$, i.e., $H = {\bf H}(F)$.

There is a maximal split torus $\bf S$ of $\bf G$, whose centralizer in $\bf G$ is a maximal torus $\bf T$ such that ${\bf B} = {\bf T}{\bf U}$, where $\bf U$ is the unipotent radical of $\bf B$. We let $\Phi$ denote the roots of $\bf G$, $\Phi^+$ the positive roots, while we also know that the choice of borel $\bf B$ fixes a set of simple roots $\Delta$. Note that a maximal parabolic ${\bf P} = {\bf M}{\bf N}$, corresponds to a subset $\theta \subset \Delta$, and $\bf N$ is generated by the root groups ${\bf N}_\alpha$ for $\alpha$ in the set of positive roots $\Phi_{\theta}^+$ of $\bf P$; we let $\rho_P$ denote half the sum of the roots in $\Phi_{\theta}^+$. We write $W(G,T)$ for the Weyl group of $\bf G$, and similarly write $W(M,T)$ for the Weyl group of $M$. Then we take
\begin{equation*}
   w_0 = w_l w_{l,M},
\end{equation*}
where $w_l$ is the longest Weyl group element of $W(G,T)$, while $w_{l,M}$ is the longest of $W(M,T)$. Furthermore, since $W(G,T) = N_G(T)/T$, we choose a representative $\tilde{w}_0 \in N_G(T)$ for $w_0$, see Section 2.1 of \cite{LomLS}. Given a representation $\rho$ of $G$, then $\tilde{w}_0(\rho)$ is the representation of $G$ defined by $\tilde{w}_0(\rho)(g) = \rho(\tilde{w}_0g\tilde{w}_0^{-1})$, $g \in G$.

We define
\[ a^*_{\bf M}=X^*({\bf M})\otimes \mathbb R, \]
where $X^*({\bf M})$ is the group of algebraic characters of ${\bf M}$. The inclusion between the corresponding split centers $A_{\bf G}\subset A_{\bf M}$ and ${\bf M}\subset {\bf G}$ define, via the restriction maps, the following decomposition
\[a^*_{\bf M}=(a^{\bf G}_{\bf M})^*\oplus a_{\bf G}^{*}\]
where $(a^{\bf G}_{\bf M})^*$ is the group of characters that are trivial on $A_{\bf G}$. Note that in our case, when $\bf M$ is a maximal Levi, the space $(a^{\bf G}_{\bf M})^*$ has dimension one. We also work with the corresponding complexifications
\[ a_{{\bf G},\mathbb{C}}^* = a_{{\bf G}}^* \otimes \mathbb{C}, \quad a_{{\bf M},\mathbb{C}}^* = a_{{\bf M}}^* \otimes \mathbb{C}, \quad (a^{\bf G}_{\bf M})_{\mathbb{C}}^*= (a^{\bf G}_{\bf M})^* \otimes \mathbb{C}. \]

We have the following map $H_M\colon M\to {\rm Hom}(a_{\bf M}^*,\mathbb R)$, defined by
\[g\mapsto(\chi\mapsto |\chi(g)|_F). \]
The kernel of $H_{M}$ is denoted by $M^1$.

Given an irreducible generic representation $\pi$ of $M = \operatorname{\textbf M}(F)$ and $\nu\in a^*_{{\bf M},\mathbb{C}}$, we define the representation $\pi_\nu$ of $M$ such that
\[ \pi_\nu(g)=q^{-\langle\nu, H(g) \rangle}\pi(g), \quad g\in M, \]
where $\langle\cdot,\cdot\rangle$ denotes the pairing between $a_{\bf M}^*$ and  ${\rm Hom}(a_{\bf M}^*,\mathbb R)$ induced by the evaluation map, and extended to $\nu\in a^*_{{\bf M},\mathbb{C}}$. 

We define $a_{\bf T}$ to be  the  dual $\mathbb{R}$-vector space of $a^*_{\bf T}=X^*({\bf T}_{\overline{F}})\otimes \mathbb R$, where $\overline{F}$ is a fixed separable algebraic closure of $F$. This space has a natural  galois action and a natural projection $p$ to $a_{\bf S}$, given by averaging the orbits. More concretely, we identify $a_{\bf S}$ with the fixed elements by the galois action of $a_{\bf T}$ and, after factoring the action through the galois group $\operatorname{Gal}(F'/F)$ where $F'$ is finite galois extension of $F$ of order $n$, the projection is then defined by the following equation
\[
n\,p(v)=\sum_{g\in \operatorname{Gal}(F'/F)}g\cdot v, \quad v\in a_{\bf T}.
\] 

Let $\alpha$ be the simple root in $\Phi$ such that $\bf P$ corresponds to the subset $\Delta\setminus \{\alpha\}$. Consider the restriction map from $a^*_{\bf T}$ to $a^*_{\bf S}$ and $\underline{\alpha}$ a simple root restricting to $\alpha$ of the root system associated with the base change of $\bf G$ to $\overline{F}$ (sometimes called the absolute root system of $G$). Moreover, $\underline{\alpha}^\vee \in a_{\bf T}$ denotes the coroot associated to $\underline{\alpha}$. Finally, we define 
\begin{equation}\label{def:alphatilde}
\widetilde{\alpha} = \dfrac{\rho_P}{\left\langle \rho_P, p(\underline{\alpha}^\vee) \right\rangle}\in (a^{\bf G}_{\bf M})^* .
\end{equation}

We use the following notation for induced representations
\[{\rm I}(\nu,\pi)={\rm Ind}_P^G(\pi_\nu), \]
where we use unitary parabolic induction. We then have a local intertwining operator on induced representations
\begin{equation*}
   {\rm A}(\nu,\pi,\tilde{w}_0) : {\rm I}(\nu,\pi) \rightarrow {\rm I}(\tilde{w}_0(\nu),\tilde{w}_0(\pi)),
\end{equation*}
which is defined via the principal value integral
\begin{equation*}
   {\rm A}(\nu,\pi,\tilde{w}_0)f(g) = \int_{N_{w_0}} f(\tilde{w}_0^{-1}ng) \, dn, \quad N_{w_0} = U \cap w_0 N_\theta^- w_0^{-1}.
\end{equation*}
It is via $\widetilde{\alpha}$ of \eqref{def:alphatilde} that we obtain an intertwining operator of a complex variable $s \in \mathbb{C}$, namely
\[
{\rm A}(s,\pi) = {\rm A}(s\widetilde{\alpha},\pi,\tilde{w}_0).
\]
Finally, in this article we only work with ${\rm N}(\nu,\pi)$ for $\nu$ ranging through the one dimensional complex space $(a_{\bf M}^{\bf G})^*_{\mathbb C}$. More precisely, if $\nu=s\tilde{\alpha}$, then
\[ {\rm N}(\nu,\pi) = {\rm N}(s,\pi), \]
where ${\rm N}(s,\pi)$ is defined via \eqref{eq:N:def}.

Central for the study of this normalized intertwining operator is a pair of local results known as the tempered $L$-function conjecture and the standard module conjecture. Both are now theorems, which we now describe.

The tempered $L$-function conjecture is about the holomorphicity of LS $L$-functions associated to tempered generic representations. It is in the setting of the Langlands-Shahidi method, where we are given a pair of quasi-split connected reductive groups $(\bf G, \bf M)$ defined over $F$, such that $\bf M$ is the Levi component of a proper maximal parabolic subgroup ${\bf P} = {\bf M}{\bf N}$ of $\bf G$. On the $L$-group side, ${}^LG$ has a corresponding maximal parabolic subgroup ${}^LP = {}^LM{}^LN$, and we let ${}^L\mathfrak{n}$ denote the Lie algebra of ${}^L\!N$. We say that $r$ is an LS-representation if it is an irreducible constituent of the representation of ${}^LM$ on ${}^L\mathfrak{n}$ obtained by conjugation. 

\subsection{Tempered $L$-functions.}\label{ss:temp:L} \emph{Suppose that $\pi$ is a generic tempered representation and $r$ is an LS representation, then $L(s,\pi,r)$ is holomorphic for ${\rm Re}(s)>0$.}

This conjecture was proved by Heiermann and Opdam when $F$ has characteristic zero in \cite{HeOp2013}, and more recently in positive characteristic by del Castillo, Henniart and  Lomel\'i  \cite{dCHeLo}.

The standard module conjecture is about the Langlands quotient for generic representations. In this context, we work with a quasi-split group $\bf G$, but $\bf M$ need not be a maximal Levi subgroup. Langlands classification, available over any local field \cite{BoWa2000}, tells us that a generic irreducible representation $\pi$ of $G$ can be expressed as the unique quotient ${\rm J}(\chi,\sigma)$ of an induced representation
\[ {\rm I}(\chi,\sigma) = {\rm Ind}_P^G(\chi\pi),
\]
where ${\bf P} = {\bf M}{\bf N}$ is a parabolic subgroup of $\bf G$ with Levi $\bf M$, and $\chi$ is a character of $M$ that is strictly positive with respect to $N$.

\subsection{Standard Module.}\label{ss:SM} \emph{If $\pi \simeq {\rm J}(\chi,\sigma)$ is a generic representation of $G$, expressed as a Langlands quotient, then 
\[ {\rm J}(\chi,\sigma) = {\rm I}(\chi,\sigma), \]
and, in particular, ${\rm I}(\chi,\sigma)$ is irreducible.}

This conjecture was proved in general over non-archimedean fields of characteric zero in Heiermann and Mui\'c \cite{HeMu2007}, and in positive characteristic by del Castillo, Henniart and Lomel\'i   \cite{dCHeLo}.

\section{Groups of classical kind and the local Ramanujan bound }\label{ss:class:Ram}

In this section we work in slightly more generality than the classical groups, and explore a bound on the Langlands parameters appearing on standard modules that we refer to as the local \emph{Ramanujan bound}. If the local representation is a component of a globally generic cuspidal automorphic representation, we prove that it satisfies the Ramanujan bound, using the approach of Kim, Zhang and the LS method in positive characteristic. 

\subsection{}\label{ss:classi:kind:def}

We say that a quasi-split connected reductive group ${\bf G}_{\ell}$ of rank $\ell$, is of \emph{classical kind} if it belongs to one of the families A, B, C, D and has the property that its maximal Levi subgroups are of the form ${\bf M} \cong {\rm Res\, GL}_m \times {\bf G}_n$, $\ell = m+n$, where ${\bf G}_n$ is a group of classical kind of rank $n$, $0 \leq n < \ell$ and ${\rm Res\, GL}_m$ is the restriction of scalars of ${\rm GL}_m$ of a quadratic extension of $F$ or is simply ${\rm GL}_m$.

The list of examples includes: split classical groups, non-split quasi-split special orthogonal and unitary groups, in addition to generalized spinor groups. We observe that it is a recursive definition, and the group ${\bf G}_n$ appearing in a maximal Levi can very well belong to a different family than the ambient group ${\bf G}_\ell$, for low rank ${\bf G}_n$. For example, the symplectic group ${\rm Sp}_{2\ell}$, has ${\rm GL}_m \times {\rm SL}_2$ as a maximal Levi, with $m = \ell-1$. Furthermore, in the above definition, we allow the case of $n=0$, where ${\bf M} \cong {\rm Res\, GL}_m$ is the Siegel Levi subgroup of ${\rm Res\, GL}_m$.

In order to be more precise, we fix a ground field $F$. We clarify that restriction of scalars is with respect to a finite separable extension $E$ of $F$, which satisfies
\[ {\rm Res}_{E/F}{\rm GL}_n(F) \simeq {\rm GL}_n(E). \]
In fact, we consider only two cases: $E=F$ or a separable quadratic extension $E/F$. For the quasi-split classical groups, in all but one case we have that $E=F$; separable quadratic extensions arising with quasi-split unitary groups.

\subsection{}\label{ss:R:bound:statement}

Given a group of classical kind ${\bf G}_\ell$ that is defined over a non-archimedean local field $F$, let $\pi$ be an irreducible generic representation of $G_\ell={\bf G}_\ell(F)$. Notice that its Levi subgroups are of the form
\[ {\rm Res}\,{\rm GL}_{n_1} \times \cdots \times {\rm Res}\,{\rm GL}_{n_b} \times {\bf G}_{n_0}. \]
Furthermore, we can fix the standard module for $\pi$ in such a way that it has the following form:
\begin{equation}\label{eq:classi:kind:std:module}
{\rm Ind} \left( \pi_b \left| \det \right|^{r_b}\otimes\cdots \otimes \pi_1 \left| \det \right|^{r_1} \otimes \pi_0 \right),
\end{equation}
where the $\pi_i$'s are tempered representations of ${\rm GL}_{n_i}(E)$, $\pi_0$ is a tempered representation of $G_{n_0}$, with ${\bf G}_{n_0}$ again a group of classical kind, and the $r_i$'s are positive real numbers, then the Langlands parameters satisfy
\begin{equation}
 0 < r_1 < \cdots < r_b. \label{eq:Rbound}  
\end{equation}
We say that a generic (unitary) representation $\pi$ of a classical kind group ${\bf G}_n$ satisfies the \emph{local Ramanujan bound} if
\begin{equation}
   r_b < 1. \label{eq:Rbound}  
\end{equation}

 In characteristic zero, Kim obtained a version of \eqref{eq:Rbound} for ${\rm SO}_{2n+1}$ and ${\rm Sp}_{2n}$, using a local-global result \cite[Lemma 3.3]{Kim2000}. For split classical groups the local Ramanujan bound is a consequence of the Lapid-Mui\'c-Tadi\'c classification of generic unitary representations of the classical groups \cite{LaMuTa2004}.
In \S~5 of \cite{Lom2019}, the local Ramanujan bound is transferred to characteristic $p$ à la Kazhdan in the case of split classical groups.

In what follows, we prove the local Ramanujan bound for groups of classical kind, under a local-global restriction. For this, we follow the aforementioned Kim's proof, but with a softer use of the local-global method, i.e. using the globalization of local supercuspidal representation of general linear groups, rather than the discrete series ones, used in Kim's original argument.

\subsection{}\label{ss:three:cases}

The setting is that of the LS machinery in our context of an an ambient group of classical kind ${\bf G} = {\bf G}_n$ of rank $n$. In the course of our proofs, there are three main cases that arise:
 \begin{enumerate}
    \item $G = {\rm  GL}_{h+k}(E)$, $M = {\rm GL}_k(E) \times {\rm GL}_h(E)$.
    \item $G = {\bf G}_{h+k}(F)$, $M = {\rm GL}_k(E) \times {\bf G}_h(F)$.
    \item $G = {\bf G}_k(F)$, $M = {\rm GL}_k(E)$.
 \end{enumerate}

Let $\rho_h$ denote the standard representation of ${\rm GL}_h(\mathbb{C})$ or one of its subgroups of classical kind. In situation (1) there is only one LS representation $r=\rho_k \otimes \tilde{\rho}_l$, hence we obtain Rankin-Selberg products by studying representations $\tau \otimes \tilde{\pi}$ of $M$ obtained from generic representations $\tau$ of ${\rm GL}_k(E)$ and $\pi$ of ${\rm GL}_l(E)$. Furthermore, we have that an LS $L$-function is equal to a Rankin-Selberg $L$-function; $L(s,\tau\otimes\tilde\pi,r)=L(s,\tau\times\pi)$ \cite{HeLo2013}.

In situation (2) we have two LS representations, and we study $\tau \otimes \tilde{\pi}$ of $M$ stemming from generic representations $\tau$ of ${\rm GL}_k(E)$ and $\pi$ of ${\bf G}_l(F)$; the first LS representation is $r={\rm std} \otimes \tilde{\rho}_l$ and the second is obtained in the next situation. Still in this second case, we denote the LS $L$-function $L(s,\tau\otimes\tilde\pi,r)$ by $L(s,\tau\times \pi).$
In (3) we again have only one LS representation arising in connection with the generic representation $\tau$ of the Siegel Levi $M={\rm GL}_k(E)$ within the ambient classical group ${\bf G}_k$; for the split classical groups we get either $r={\rm Sym}^2 \rho_k$ or $\wedge^2 \rho_k$, cf. \cite{HeLo2011}. 

\subsection{}\label{ss:Kim:local-global} We are now ready to prove our local-global result, where ${\bf G} = {\bf G}_n$ is a group of classical kind defined over a global function field $K$; hence, $\bf G$ is also defined over each local field $F$ obtained from $K$ after completion at a place $v$, $F = K_v$, in addition to the ring of ad\`eles $\mathbb{A}_K$ of $K$.

\noindent{\bf Proposition.} \emph{Let $\pi$ be a generic representation of ${\bf G}_n(F)$ that arises as the local component of a globally generic cuspidal automorphic representation $\pi' = \otimes_v \pi_v'$of ${\bf G}_n(\mathbb{A}_K)$. More specifically, there exists a place $v_0$ of a global function field $K$, such that $F = K_{v_0}$ and $\pi_{v_0}' = \pi'$. Then, $\pi$ satisfies the local Ramanujan bound, i.e. $r_b$ arising from the representation $\pi$ of ${\bf G}_{n}(F)$ as in \eqref{eq:classi:kind:std:module} satisfies}
\[ r_b < 1. \]
\begin{proof}First, we recall that every tempered representation of a general linear group is fully induced from discrete series representations. Let us apply this to  the standard module of $\pi$  \eqref{eq:classi:kind:std:module}. After replacing every $\pi_i$ with the corresponding induced representation from discrete series representations, we can write $\pi$ as the full induced representation of the form
\begin{equation}\label{eq:stdmodule:discrete}
{\rm Ind} \left( \delta_b \left| \det \right|^{r'_b}\otimes\cdots \otimes \delta_1 \left| \det \right|^{r'_1} \otimes \pi_0 \right).
\end{equation}
In this last equation, the $\delta_i$'s are discrete series representations of ${\rm GL}_{m_i}(E)$, $\pi_0$ is a tempered representation of $G_{n_0}$ for a group of classical kind ${\bf G}_{n_0}$, and the Langlands parameters satisfy $r_b'=r_b, r'_1=r_1$ and
\begin{equation*}
 0 < r'_1 \leq \cdots \leq r'_b.   
\end{equation*}
Now, we look into the possibilities for the discrete series appearing in \eqref{eq:stdmodule:discrete}, focusing on $\delta_b$. We look at supercuspidal representations arising in connection to Bernstein-Zelevinsky classification of quasi-discrete series representations \cite[Theorem 9.3]{Ze1980}. Indeed, the discrete series $\delta_b$ is the Bernstein-Zelevinsky representation $L(\Delta)$ associated to a segment of supercuspidal representation $\rho$:
\[\Delta=\{\rho,\nu\rho,\dots,\nu^{u-1}\rho\}, \quad u\in \mathbb Z_{>0}.\]

In order to proceed further, we study their  associated $L$-functions. The properties of Rankin-Selberg $L$-functions that we need can be found in \cite{JaPSSh1983}, however, we use \cite{He2002}, since it better suits our notations. First, we claim that
\[ L(s,\rho\times  \Tilde \delta_b)=L(s,\rho\times\Tilde \rho). \]
We recall that $\Tilde \delta_b=\nu^{1-u}L(\Tilde\rho,\nu\Tilde\rho,\dots,\nu^{u-1}\Tilde\rho)$. Now, we can compute the $L$-function in terms of the Bernstein-Zelevinsky classification \cite[Section 2.6]{He2002}, as follows
\begin{align*}
L(s,\rho\times  \Tilde\delta_b) &= L(s+1-u, \rho \times L(\Tilde\rho,\nu\Tilde\rho,\dots,\nu^{u-1}\Tilde\rho)) \\
&= L(s+1-u,\nu^{u-1}\rho\times\Tilde \rho)
L(s,\rho\times\Tilde \rho).
\end{align*} 
Furthermore, we claim that  $L(s-r_b', \rho\times\Tilde\rho)$ is a factor of $L(s,\rho\times \pi )$. Indeed, using the multiplicativity property of LS $L$-funtions in positive characterisitc \cite{LomLS}, and the discussion above, we have that,
\begin{align*}
L(s,\rho\times\pi) = L(s-&r_b',\rho\times \Tilde\delta_b)\cdots L(s-r_1', \rho\times \Tilde\delta_1) \\
& \times L(s,\rho\times\pi_0) L(s+r_1',\rho\times\delta_1)\cdots L(s+r_b',\rho\times\delta_b).
\end{align*}
Hence, combining the previous equality, we get that $L(s-r'_b, \rho\times\Tilde\rho)$ is also a factor of $L(s,\rho\times\pi)$.

On the other hand, there exists a globally generic cuspidal $\rho' = \otimes_v \rho_v'$ of a general linear group such that $\rho'_{v_0}\cong \rho$ and $\rho'\otimes\pi'\not \cong w_0(\rho'\otimes\pi')$. In this situation, we observe that Proposition 5.4.2 of \cite{dC} applies to $\rho'\otimes\pi'$, in the context of groups of classical kind, to get that $L(s,\rho\times\pi)$ is holomorphic for $\operatorname{Re}(s)\geq 1$. Therefore, since the factor $L(s-r_b', \rho\times\Tilde\rho)$ of $L(s,\rho\times\pi)$  has a pole at $s=r_b'=r_b$ \cite[Section 3.3]{He2002}, we must have that $r_b<1$. 
\end{proof}

\subsection*{Remark} The above proposition is valid in the situation where there exists a place $v_0$ of a global function field $K$, together with an isomorphism $F \simeq K_{v_0}$ such that the component $\pi_{v_0}'$ of $\pi'$ is isomorphic to $\pi$.

\section{Proof of Kim's Assumption A}\label{sec:kimproof}

We prove our main result, Kim's Asssumption A, for groups of classical kind in characteristic $p$
under the condition that the local Ramanujan bound $r_b < 1$ is satisfied.  As a consequence, we establish the result for split classical groups, where the Ramanujan bound holds true, and in general under the local-global condition of the previous section.

\subsection{}
Let us first recall the validity of Kim's Assumption A for general linear groups, where we have a result of M\oe glin and Waldspurger. More precisely, Proposition I.10 of \cite{MoWa1989} tells us that given tempered representations $\tau$ of ${\rm GL}_n(E)$ and $\tau'$ of ${\rm GL}_{n'}(E)$, then ${\rm N}(s,\tau \otimes \Tilde{\tau}')$ is holomorphic and non-zero on ${\rm Re}(s) > -1$. Note that we take the contragredient of $\tau'$ to match with the LS factors with the Rankin-Selberg ones \cite{HeLo2013}. Then Kim's local Assumption A follows by the classification of unitary generic representations for general linear groups \cite{Ta1986}. This is a key ingredient for the more general cases.

We now follow \cite{Kim1999, Kim2000, Kim2005,Zh} in order to prove Kim's Assumption A for groups of classical kind; under the local-global restriction in general, and with no restriction for the split classical groups.

\noindent{\bf Theorem.} \emph{Let $\pi$ be a generic representation of ${\bf G}_n(F)$. Suppose that $\pi$ satisfies the local Ramanujan bound. Then ${\rm N}(s,\tau \otimes \tilde{\pi})$ is holomorphic and non-zero on ${\rm Re}(s) \geq 1/2$ for every generic unitary representation $\tau$ of ${\rm GL}_m(E)$.}

\begin{proof}
We first recall the results of \S~7 of \cite{Ta1986}, stating that a generic unitary representation $\tau$ of ${\rm GL}_m(E)$ can be written as the full induced representation
\begin{equation}\label{eq:GL:std:module}
    {\rm Ind}_Q^{\operatorname{GL}_m(E)}(\xi_1|\det|^{t_1}\otimes \cdots \otimes \xi_d|\det|^{t_d}\otimes \xi_{d+1} \otimes\xi_d|\det|^{-t_d} \otimes \cdots \otimes\xi_1|\det|^{-t_1}),
\end{equation}
where the $\xi_i$'s are tempered representations of ${\rm GL}_{m_i}(E)$ and
\[ 0<t_1\leq \cdots \leq t_d<1/2. \] 
Note that $d$ could be zero and that  $\xi_{d+1}$ could not appear, in which case $\tau$ was already a tempered representation. The representation $\pi$ has a standard module, given by \eqref{eq:classi:kind:std:module}.

We now observe that the multiplicativity property of LS local factors, can be used to obtain the corresponding property for normalized operators. This leads to the following formula:
\begin{align*}
    {\rm N}(s,\tau\otimes\Tilde \pi)=& \prod_{j=1}^{b}\prod_{i=1}^d{\rm N}(s,\xi_i \otimes\tilde{\pi}_j|\det|^{t_i-
    r_j}) {\rm N}(s,\xi_i \otimes\tilde{\pi}_j|\det|^{- t_i-
    r_j})  \\
    &  \prod_{j=1}^{b}\prod_{i=1}^d{\rm N}(s,\xi_i \otimes{\pi_j}|\det|^{ t_i+
    r_j}) {\rm N}(s,\xi_i \otimes{\pi_j}|\det|^{- t_i+
    r_j}) \\
    &  \prod_{j=1}^{b}{\rm N}(s,\xi_{d+1}\otimes{\tilde \pi_j}|\det|^{ 
    -r_j}){\rm N}(s,\xi_{d+1}\otimes{\pi_j}|\det|^{
    r_j})\\
    &{\rm N}(s,\xi_{d+1} \otimes\tilde{\pi}_0)\prod_{i=1}^d {\rm N}(s,\xi_i \otimes\tilde{\pi}_0|\det|^{ t_i}) {\rm N}(s,\xi_i \otimes\tilde{\pi}_0).
\end{align*}
Here, we observe that each factor appearing on the right hand side falls into one of the three Langlands-Shahidi situations for pairs $({\bf G},{\bf M})$, discussed in \S~\ref{ss:three:cases}. In fact, only the first two cases arise.

In the first situation, there are two representations of general linear groups. If $\operatorname{Re}(s) \geq 1/2$, in all combinations for $i$ and $j$, we see that $\operatorname{Re}(s\pm t_i\pm r_j )>-1$. Then the holomorphy and non-vanishing of the normalized intertwining operator corresponding to $\xi_i|\det|^{\pm t_i} \otimes \Tilde \pi_j|\det|^{\pm r_j}$ holds thanks to Proposition I.10 of \cite{MoWa1989}. 

In the second situation of a $\pi_0$ and $\xi_i$. If $\operatorname{Re}(s)\geq 1/2$, then $\operatorname{Re}(s\pm t_d)\geq 0$. Since $\pi_0$ and $\xi_i$ are tempered (or $\tau$ is tempered, i.e. in the case when $d=0$ and $\xi_{d+1}$ does not appear), we get that it is holomorphic for $\operatorname{Re}(s)\geq 0$  thanks to \ref{ss:temp:L} Theorem. Indeed for $\operatorname{Re}(s)>0$, the (non-normalized) intertwining operator for tempered representations is holomorphic and non-vanishing for $\operatorname{Re}(s)>0$ \cite[Equation (10) and Proposition IV.2.1]{Wa2003}. The holomorphicity for $\operatorname{Re}(s)=0$, we proceed exactly as in \cite[Lemma 2]{Zh}, as all the steps there are available in characteristic $p$. 

Finally, for the non-vanishing property at $\operatorname{Re}(s)=0$ in the second situation, we are in the special case of a maximal parabolic subgroup in Theorem~3 of \cite{Zh}. Zhang's argument is general, making use of combinatorial and basic geometric properties or root spaces, together with the multiplicativity property provided by the Langlands-Shahidi method, now available in positive characteristic, to deduce non-vanishing. In our notation, we have that  the normalized intertwining operator satisfies the following formula
\[
{\rm Id}={\rm N}(\Tilde w_0(s\Tilde{\alpha}),\Tilde w_0(\xi_i\otimes\tilde{\pi}_0)){\rm N}(s,\xi_i\otimes\tilde{\pi}_0), \]
where $\Tilde w_0$ is the element of \S~\ref{s:prelim}. More precisely, this corresponds to the case $\theta'=\Delta$ exactly as in \cite[Theorem 3, p.393]{Zh}. With this formula, the non-vanishing of ${\rm N}(s,\xi_i\otimes\tilde\pi_0)$ for $\operatorname{Re}(s)=0$, follows form the holomorphichty at ${\rm Re}(s)=0$ of both terms on the right hand side.
\end{proof}

\subsection{} In the course of the above proof, we can observe that only the cases (1) and (2), introduced in \S~\ref{ss:three:cases} for groups of classical kind, are used. However, let us now prove the corresponding result for case (3).

We first reduce to the supercuspidal representation case, via the Bernstein-Zelevinsky classification. The representation $\tau$ is fully induced from discrete series, as in \eqref{eq:stdmodule:discrete}. And then, we write each of those discrete series using segments of supercuspidal representations. Now, if our representation is unitary supercuspidal and using that the Intertwining operator has at most simple poles at $\operatorname{Re}(s)=0$ \cite[Proposition IV.1.2]{Wa2003}, we get that the Normalized Intertwining operator is holomorphic and non-zero outside of $\operatorname{Re}(s)=-1/2,-1$, using the same arguments in \cite[Lemma 4.1]{Kim2005}. We conclude by following the exact arguments as in \cite[Lemma 3.3, Proposition 3.4]{Kim1999}. We summarize what we have shown in the setting of case (3):

\begin{quote}
\emph{The normalized intertwining operator ${\rm N}(s,\tau)$ is holomorphic and non-zero on ${\rm Re}(s) \geq 1/2$ for every generic unitary representation $\tau$ of ${\rm GL}_m(E)$.}
\end{quote}

\medskip

\noindent{\sc \Small H\'ector de Castillo, Departamento de Matem\'atica y Ciencia de la Computaci\'on, Universidad de Santiago de Chile, Las Sophoras 173, Estación Central, Santiago}\\
\emph{\Small E-mail address: }\texttt{\Small hector.delcastillo@usach.cl}

\noindent{\sc \Small Luis Alberto Lomel\'i, Instituto de Matem\'aticas, Pontificia Universidad Cat\'olica de Valpara\'iso, Blanco Viel 596, Cerro Bar\'on, Valpara\'iso, Chile}\\
\emph{\Small E-mail address: }\texttt{\Small Luis.Lomeli@pucv.cl}

\end{document}